# Differential and Integral Calculus for Logical Operations
# A Matrix-Vector Approach

Eduardo Mizraji
Group of Cognitive Systems Modeling, Biophysics Section,
Facultad de Ciencias, Universidad de la República
Montevideo, Uruguay

Address:
Dr. Eduardo Mizraji
Sección Biofísica, Facultad de Ciencias, UdelaR
Iguá 4225, Montevideo 11400, Uruguay
e-mails: 1) emizraji@gmail.com,  2) mizraj@fcien.edu.uy
Phone-Fax: +598 25258629

ABSTRACT

A variety of problems emerged investigating electronic circuits, computer devices and cellular automata motivated a number of attempts to create a differential and integral calculus for Boolean functions. In the present article, we extend this kind of calculus in order to include the semantic of classical logical operations. We show that this extension to logics is strongly helped if we submerge the elementary logical calculus in a matrix-vector formalism that naturally includes a kind of fuzzy-logic. In this way, guided by the laws of matrix algebra, we can construct compact representations for the derivatives and the integrals of logical functions. Inside this semantic-algebraic calculus, we obtain expressions for the derivatives of some of the basic logical operations and show the general way to obtain the derivatives of any well-formed formula of propositional calculus. We show that some of the basic tautologies (*Excluded middle, Modus ponens, Hypothetical syllogism*) are members of a kind of hierarchical system linked by the differentiation algorithm. In addition using the logical derivatives we show that relatively complex formulas can collapse in simple expressions that reveal clearly their hidden logical meaning. The search for the antiderivatives produces naturally an integral calculus. Within this logical formalism an indefinite integral can always be found for any logical expression. Moreover, particular integrals can be constructed based on detachment properties that lead to logical expressions of growing complexity. We show that these particular integrals have some similarities with the "generalizing deduction" procedures investigated by Łukasiewicz.

**Keywords**: Boolean derivative, Logical antiderivative, Matrix-vector logics, Many-valued logic

1. INTRODUCTION.

In the time of Leibniz and Newton, the development of the infinitesimal calculus was strongly stimulated by mathematical and physical problems that belonged, mainly, to the field of continuous mathematics. A great amount of our present scientific knowledge is a consequence of this development. Yet, in our times many problems of theoretical and technical importance belong to the territory of discrete mathematics. We mention as examples the optimization of logical circuits or the design of cellular automata. In these situations, the direct transposition of the methods of infinitesimal calculus is not possible. For this kind of problems an operation called Boolean derivative has been defined [8, 9, 12, 18, 24, 42, 43]. The Boolean derivative of a formula with respect to a variable gives the condition under which the formula changes its value whenever the variable does. This derivative produces a new Boolean function that only depends on the remaining variables [42, 43]. After the proposal by Shannon in his MSc thesis [38], the theory and analysis of electrical circuits has been largely improved using the formalism of binary variables and the use of what is usually called "Boolean logic" (even if some of the main operations of this binary logic contradict the definitions of Boole [3, 14]).

Let us mention that the debates about the nature of Boolean operations were extremely fertile and produced important theoretical and practical ramifications. Symmetric difference was partially implicit in the work of Boole (see, for instance, [3], pp. 32-33, Dover Edition); in fact, this operation was adopted by a number of mathematicians in their formalizations of Boolean logical theory. For instance, inspired by Boole's theory, M. H. Stone published in 1935 [40] a Boolean algebra based on two fundamental operations: the multiplication and the symmetric difference. Inclusive addition is defined from these two basic operations. These operations allow him to show that his Boolean algebra was a particular class of ring. The use of the symmetric difference (exclusive-or) in Boolean representations was central for the construction of a Boolean differential calculus.

The powerful instrument of binary representations of electrical circuits and the problems of predicting the effects of relays and switches, naturally leaded to the definition of a discrete derivative. This "Boolean derivative" allows detection of the effect of variable switching on the behaviour of the circuit [8, 9, 12, 24, 42]. A complementary creation is the Boolean integral, also used in relation with problems concerning electrical circuit theory [9, 42]. Some of the approaches looking for Boolean differential and integral calculus have used the formalism of matrix algebra (see, in particular, the recent work by Cheng [9]). Cheng adopts the matrix algebra to represent a consistent logical formalism, assuming that the variables participating in it are vectors and matrix built up using binary numbers.

An important domain of application of the Boolean derivatives is the theory of cellular automata as developed by Vichniac [43]. Even for the simplest situations, the prediction of the dynamic behavior of these mathematical devices is a hard problem. In particular, the Boolean derivative has been used to analyze the structure of the space of 256 logical functions that defines the two-states, two-neighbour 1-dimensional cellular automata (Wolfram automata of the class k = 2, r = 1) [33, 43]. Conversely, very complex problems of continuous mathematics (the solutions of Navier-Stokes equations) can be extremely well approximated using cellular automata governed by discontinuous laws [45]. This fact suggests the existence of links, not yet well understood, between continuous and discrete representations. Hence, the development of a differential and integral calculus defined on logical variables finds a first motivation in the practical problems of discrete mathematics, and in their connection with non-linear continuous models.

Another reason to explore the properties of a differential and integral calculus of discrete mathematical structures is related to the fact that in classical analysis differentiation and integration show natural links with the losing or gaining of information. In the case of real functions, the differentiation represents an irreversible process: it is executed by a well defined algorithm, and produces a loss of information in the sense that differentiation is a many-to-one application. Consequently, the integration, the reverse operation, is not single-valued, and we do not have general algorithms to obtain antiderivatives for real-valued functions [16, 17]. In the simplest cases, the construction of indefinite integrals uses a basic set of elementary antiderivatives, plus theorems that govern their transformations. Moreover, in some cases the antiderivatives exist, but they are not expressible in terms of finite representations as is well illustrated by Gaussian exponential function $\exp(-x^2)$. The fact that a consistent definition for the Boolean derivatives and antiderivatives has been proposed by many authors, opens the possibility to explore these topics in the discrete domain.

In the theories previously referred, Boolean derivatives and integrals are applied to canonical Boolean expressions based on additions, multiplications and complementarities. Clearly, these Boolean expressions are always behind any well-formed formula (wff) of propositional calculus, including the basic logical operators. In a sense we can say that canonical Boolean expressions are a syntax that supports the constructions of the logical "semantics" implied in any wff of propositional calculus.

Our purpose here is to adapt the Boolean differential and integral calculus to the semantic expressions of propositional logic. For example, if we have a standard formulation of the classical tautology *modus ponens*, $\left[ p \wedge (p \rightarrow q) \right] \rightarrow q$, we can obtain a derivative able to be expressed as follows:

$$\frac{\partial}{\partial p}\left\{\left[p \wedge (p \rightarrow q)\right] \rightarrow q\right\} \equiv \neg(q \vee \neg q).$$

To emphasize the interest of conserving the semantics of logical operations, we observe that this derivative of the *modus ponens* provides the negation of other basic tautology: the *excluded middle*.

In the case of integration, we show later in this article that for any wff we can define a set of general integrals and particular integrals, all them reciprocal operation of the Boolean differentiation. For instance, given the logical formula $p \rightarrow \neg q$, we can show that a general integral for this expression is given by

$$\int (p \rightarrow \neg q) \partial r \equiv \left[ (p \rightarrow \neg q) \rightarrow r \right],$$

and a particular integral for the same expression is given by

$$\int_P (p \rightarrow \neg q) \partial r \equiv \left[ (p \vee r) \rightarrow \neg (q \vee r) \right].$$

The classical methods of binary logical can always be used to obtain and verify these expressions, Nevertheless, the processing of differentiation and integration over logical expressions can be carried out with a different methodology if we submerge the logical formalism into the realm of linear algebra. The procedure employed here is an example of the "bypasses" analyzed by Z.A. Melzak [27]: the problem to differentiate or integrate any wff logical formula is mapped into the domain of matrix algebra; in this domain the definitions are established and their consequences are investigated; finally, as the third step of Melzak's bypass, the results can be re-translated into the classical format of propositional logic. It has been previously shown that the matrix-vector representation of propositional logic (formalism named "vector logic" [30]), allows to perform the logical calculus as operations of linear algebra. In this formalism the basic logical functions (as negation, disjunction, implication, etc.) are represented by matrices and the truth-values are represented by vectors. Within this vector logic, the equivalences of propositional calculus become equalities between algebraic equations, and the logical symbols become matrix operators. One of the remarkable facts of this vector logic is that some basic logical equivalences (eg, the De Morgan's laws) become identities between matrix operators, not dependent on the vector truth-values [30, 32, 35]. The matrix-vector formalism of this vector logic has very close relations with some matrix models of neural associative memories, and in this sense they establish a promising point of departure to investigate the neural bases of human reasoning [29, 36].

The matrix-vector approach has another particularity: in this algebraic calculus, uncertainties in the truth-values can be expressed as linear combinations of the basic truth-values affected by probabilistic weights. The matrix operators are composed by real numbers, are defined over two basic truth-value vectors (each vector corresponding to each one of the binary truth-values) and programmed to produce a binary vector logic. However, these "binary" matrices based on two vectors but composed by real numbers, are capable of processing the probabilistic truth-values and, in this way, to produce a particular class of many-valued logic [30, 35]. This fact allows extending the range of the logic differential and integral calculus out of the binary domain, and gives a further insight to interpret the interaction of differentiation and integration with the losing or gaining of information.

The structure of the paper is the following. In the first parts we review the basic operations of elementary propositional calculus, with a special emphasis in their relation with the Boolean polynomials. Then we show how vector logic is a matrix-vector translation of this formalism created by Boole, into the language and operations of matrix algebra. We also illustrate how this vector logic produces a probabilistic scalar many-valued logic in the presence of probabilistic inputs. We adapt to this formalism the definitions of Boolean derivatives and apply it to the basic logical operators. Cross derivatives and successive derivatives are also described. We show how successive differentiation of classical tautologies discovers the existence of an interesting semantic hierarchy between them, and how differentiation can be an useful mean to simplify arguments. Finally, looking for antiderivatives, we construct firstly a general integral only dependent on the basic function considered, and secondly a particular integral that also depends on arbitrary logical expressions.

## 2. LOGICAL FUNCTIONS
### 2.1 The classical logical functions and the Boolean Polynomials

Classical binary logic is based on a reduced number of mathematical functions depending on one (monadic) or two (dyadic) variables. In a binary base set $\{1, 0\}$, the value 1 corresponds to "true" and the value 0 to "false". The monadic functions are of the form $y = M(x)$, and the dyadic functions are of the form $z = B(x, y)$ with $x, y, z \in \{1, 0\}$. Table 1 shows two monadic functions, and Table 2 displays the most important dyadic functions.

**Table 1: Monadic functions**

| x | ID | NOT ($\neg$) |
|---|----|----|
| 1 | 1  | 0  |
| 0 | 0  | 1  |

**Table 2 : Dyadic functions**

| x | y | AND ∧ | OR ∨ | IMPL → | NAND | NOR | XOR | EQUI ≡ |
|---|---|---|---|---|---|---|---|---|
| 1 | 1 | 1 | 1 | 1 | 0 | 0 | 0 | 1 |
| 1 | 0 | 0 | 1 | 0 | 1 | 0 | 1 | 0 |
| 0 | 1 | 0 | 1 | 1 | 1 | 0 | 1 | 0 |
| 0 | 0 | 0 | 0 | 1 | 1 | 1 | 0 | 1 |

In the framework of the classical logic, the tautologies are logical expressions that produce the truth-value 1 (or "true") for all the possible values of their logical variables. Some of the basic tautologies are the following ($p,q,r \in \{0,1\}$):

*Excluded middle*: $p \vee \neg p$

*Modus ponens*: $[p \wedge (p \rightarrow q)] \rightarrow q$

*Hypothetical syllogism*: $[(p \rightarrow r) \wedge (r \rightarrow q)] \rightarrow (p \rightarrow q)$

It is remarkable that when George Boole established the development of logical operations as polynomials, he established a firm bridge between classical logics and algebra [3]. For the case of monadic operators, the Boolean polynomial looks as follows:

$$f(x) = f(1)\, x + f(0)(1-x)$$

The 4 different monadic operations result from the different binary values for the coefficients. For instance the ID operation requires $f(1) = 1$ and $f(0) = 0$, and NOT happens if $f(1) = 0$ and $f(0) = 1$. For the case of dyadic operators, the Boolean polynomials are of the general form

$$f(x,y) = f(1,1)xy + f(1,0)x(1-y) + f(0,1)(1-x)y + f(0,0)(1-x)(1-y),$$

and the different combinations of binary coefficients generate the 16 dyadic logical functions. The operations illustrated in the Table 2 can be translated to this polynomial format when the coefficients take the values indicated in the table. For instance: NAND requires that $f(1,1) = 0$ and $f(1,0) = f(0,1) = f(0,0) = 1$. These Boolean polynomials can be immediately extended to any number of variables, $f(x_1,\ldots,x_i,\ldots,x_p)$, producing a large potential variety of logical operators [35].

The complex relations between mathematics and logics, and the way many concepts change of status in different epochs, have been reviewed in [20]. This review shows the critical role of invention for producing innovative advances and sustained challenges in mathematics and logic. In the first paragraph of an influential book originally published in 1905 by Louis Couturat, we have a clear statement of a mathematician about the calculus of the logic developed by Boole: *"Les lois fondamentales de ce calcul on été inventées pour exprimer les principes du raisonnement, les 'lois de la pensée'; mais on peut considerer ce calcul au point de vu purement formel, qui est celui des Mathématiques, comme une Algèbre reposant sur certaines principes arbitrairement poses"* [11] (see footnote 1). The logical polynomials created by Boole became a subject of interest for many mathematicians and logicians. In Section 24 of his book Couturat analyze Boolean polynomials and adapt the original formalism to his pure algebraic approach. He explicitly uses the dual interpretation of the polynomial variables as categories or propositions, in both cases bounded by 0 and 1. The formalization of Boolean algebra was expanded in many directions and a large number of refinements were proposed in the next decades. After his short paper of 1935 [40] Stone published a comprehensive and highly influential theory of Boolean rings [41]. This theory is a remarkable proof of the argument given in [20.] showing how some initially "heterodox" inventions, as the idempotency of Boolean variables, become basic for a solid and fertile theory. Recent developments expand the theory of Boolean polynomials creating a polynomial ring calculus that embraces many aspects of the theory of logics (eg. representation of syllogism and many valued logics) [5]. Using this kind of polynomial rings has been possible to develop a novel semantic procedure to investigate logical modalities [2].

2.2. **The Operators of Vector Logic**.
Vector logic [30] is an algebraic model of elementary logic based on matrix algebra. It is assumed that the truth-values map on Q-dimensional vectors and that the monadic and dyadic operations are executed by matrix operators. It is a remarkable property of the vector logic formalism that when the matrices operate over vectors representing truth-values, the generated formulas show the same order of operations and logical variables produced by Polish notation in the case of non-vectorial variables [32, 35]. This fact illustrates the implicit operator approach inside Łukasiewicz conception of logical formalism [25]. In addition, our representation generates an operator theory where the logical operators are themselves subject to the laws of matrix algebra (this point has been investigated with detail in [32]). It is well known that George Boole was an expert in the application of operators in the domain of differential and difference equations [14, 4]. Since the algebraic properties of matrices were communicated by Cayley in 1858 [7], unfortunately Boole had no occasion (as far as we know) to connect the newborn matrix theory with his algebraic formalism for the logic. It is a peculiar historical fact that Boole and Cayley maintained a brief and, in a sense, slightly divergent correspondence about the algebra of logic between 1847 and 1855 (reproduced in [4]). The potentialities of matrix algebra for the theory of the logic was rapidly anticipated by Charles Peirce around 1870 and explored by Irving Copi in 1948 (references in [10]). The matrix formalism has revealed the interesting possibility of representing the deep ideas of Łukasiewicz about many-valued logics and modalities as a full operator theory [31, 35].

---

1 "The fundamental laws of this calculus were devised to express the principles of reasoning, the 'laws of thought'. But this calculus may be considered from the purely formal point of view, which is that of mathematics, as an algebra based upon certain principles arbitrarily laid down". Translated to English by Lydia Gillingham Robinson and published online by The Project Gutenberg: EBook of The Algebra of Logic, by Louis Couturat, Release Date: January 26, 2004 [EBook #10836].

A general approach to this vector logic has been described in [35]. Here, we will be mainly concerned with monadic and dyadic operators. We want to mention that different formalisms for the logic based on matrices and vectors have been investigated in relation to the basic theory of Boolean functions [8, 9]. quantum physics [13, 28] and new physical computing procedures [44]. In addition, a recent article [21] describes and enlarges the pioneering contributions of G.N Ramachandran showing the accuracy of the matrix-vector formalism to capture some aspects of Indian Logic.

A propositional calculus can be considered as a classification system that assigns a truth-value to each proposition. In the traditional binary logic, the truth-values are: true t (or "yes") aand false f (or "not"), and the basic set for the definition of mathematical logic functions is $\tau_2 = \{t, f\}$. This binary logic possesses 4 monadic operations and 16 dyadic operations. The construction of a two-dimensional vector logic begins establishing of a correspondence between the truth-values t and f, and two Q-dimensional normalized column vectors: $t \mapsto s$ and $f \mapsto n$, $s, n \in \mathbb{R}^{Q \times 1}$, with $Q \geq 2$ (the vector notation using "s" and "n" is based in the Spanish "yes", SI, and "not", NO.). This correspondence produces a set of vector truth-values:

$$V_2 = \{s, n\}.$$

The logical operations defined over this set of vectors lead to matrix operators. It is especially interesting the fact that the set of truth-values $V_2 = \{s, n\}$ and the associated matrix logical gates, allow to compute linear combinations of vectors s and n, that become a first natural representation for uncertain truth-value assignments.

To simplify the notation and the results, we are going to assume here that s and n are orthonormal column vectors, but this assumption is not necessary in general, and in [32] we show how to derive the logical operators from linear independent vectors using the Moore-Penrose pseudoinverse of a matrix. We want to mention that for the following results we only need to consider vectors of dimension 2. Nevertheless, we want to retain the general presentation based on Q-dimensional vectors because this approach keeps open the possibility to extend, in future works, many aspects of this formalism to a larger number of vector truth-values. Assuming more than two vector truth-values, allows us to explore some interesting aspects of many-valued logics: on the one hand, we can add new vector logic values, orthogonal to the two basic vector truth-values s and n; on the other hand, this enlarged set of vectors can admit inputs to the matrix operators weighted with scalar probabilistic coefficient. We begun this kind of exploration in a paper published in 2008 [35].

### 2.2.1. Basic operations

The scalar product between Q-dimensional column vectors, $u^T v = \langle u, v \rangle$, is the operation responsible of the properties displayed by vector logic. The orthonormality between vectors s and n implies that $\langle u, v \rangle = 1$ if $u = v$, and $\langle u, v \rangle = 0$ if $u \neq v$, $u, v \in V_2$.

### a) Monadic Operators

The basic monadic operators for this two-dimensional vector logic are generated by the mapping

$$\text{Mon}: V_2 \to V_2.$$

This mapping produces four square matrices $I, K, M, N \in \mathbb{R}^{Q \times Q}$. The matrices I and N are respectively the identity and the negation matrices, and K and M are two operators that produce a constant output [32].

**a1) Identity.** A logical identity $ID(p)$ produces a matrix behaving as follows: $Iu = u$, $u \in V_2$, and the structure of this matrix is

$$I = ss^T + nn^T.$$

Hence, due to the orthogonality of s respect to n, we have $Is = s\langle s, s \rangle + n\langle n, s \rangle = s$ and $In = n$.

**a2) Negation.** The classical negation $\neg p$ is represented by the matrix operation $Nu$, $u \in V_2$, with

$$N = ns^T + sn^T.$$

Consequently, $Ns = n$ and $Nn = s$. Note that the involutory behavior of the logical negation, $\neg(\neg p) \equiv p$, corresponds with the fact that $(N)^2 = I$ (the vector logic identity matrix is not generally an identity matrix in the sense of matrix algebra, except in particular cases).

a3) **Constant operators**. The following matrices give monotonic outputs:

$$K = ss^T + sn^T$$
$$M = ns^T + nn^T.$$

Consequenttly, $Ks = Kn = s$ and $Ms = Mn = n$,

b) **Dyadic operators**
The 16 two-valued dyadic operators correspond to the following mapping:

$$\text{Dyad}: V_2 \otimes V_2 \to V_2 \ .$$

This mapping generates the rectangular matrices $T \in \mathbb{R}^{Q \times Q^2}$.

The different matrices T that execute these dyadic operations are based on the properties of the Kronecker product. As we described previously (see [30] and [32]), the Kronecker product allows representing logical variables without the need of modular arithmetic. Some of the recent approaches to logic formalism that use matrices and Kronecker products adapt the modular arithmetic for the algebraic operations between matrices and vectors (see, for instance, [8, 9]). In the simplest version of the formalism used in the present paper, real vectors play the role of Boolean variables due to orthogonality; in the case of non-orthogonal linearly independent vectors, matrix pseudoinverses generate similar results [32]. A remarkable aspect of this algebraic representation is the natural emergence of a matrix version of Boolean polynomials; in these polynomials, the matrix-vector operations and the Kronecker products replace modular addition and multiplication [32, 35].

We summarize in what follows the definition and some basic properties of this product [19]. Given two matrices $A = [a_{ij}]_{m \times n}$ and $B = [b_{ij}]_{p \times q}$, the Kronecker product $A \otimes B$ is given by

$$A \otimes B = [a_{ij}B]_{(mp) \times (nq)}.$$

Two properties of this product are essential for the formalism of vector logic:
(P1) $(A \otimes B)^T = A^T \otimes B^T$
(P2) $(A \otimes B)(A' \otimes B') = (AA') \otimes (BB')$ .

Property (P2) needs conformable matrices (or vectors). For two r-dimensional column vectors a and c, and two r'-dimensional vectors b and d, (P2) implies
$(a \otimes b)^T (c \otimes d) = (a^T c)(b^T d) = \langle a, c \rangle \langle b, d \rangle$.

The matrix versions of the basic dyadic operators are described in the next paragraphs.

**b1) Conjunction.** The conjunction between two propositions $p \wedge q$ is represented by a matrix that acts on two vector truth-values: $C(u \otimes v)$ , $u, v \in V_2$ . This matrix C reproduces the features of the classical conjunction truth-table:

$$C = s(s \otimes s)^T + n(s \otimes n)^T + n(n \otimes s)^T + n(n \otimes n)^T$$

and operates as follows: $C(s \otimes s) = s$ ; $C(s \otimes n) = C(n \otimes s) = C(n \otimes n) = n$ .

**b2) Disjunction.** The classical disjunction $p \vee q$ is executed by the matrix

$$D = s(s \otimes s)^T + s(s \otimes n)^T + s(n \otimes s)^T + n(n \otimes n)^T ,$$

being $D(s \otimes s) = D(s \otimes n) = D(n \otimes s) = s$ and $D(n \otimes n) = n$ .

In the two-valued logic, the conjunction and the disjunction operations satisfy the De Morgan Law: $p \wedge q \equiv \neg(\neg p \vee \neg q)$ (and also the dual: $p \vee q \equiv \neg(\neg p \wedge \neg q)$ ). For two-dimensional vector logic this Law is also verified

$$C(u \otimes v) = ND(Nu \otimes Nv) .$$

The Kronecker product allows the following factorization:

$$C(u \otimes v) = ND(N \otimes N)(u \otimes v).$$

A remarkable fact is the following. We can prove directly, from the previous matrix definitions, that in the two–dimensional vector logic the De Morgan Law is a law involving operators, and not only a law concerning operations:

$$C = ND(N \otimes N) .$$

For a detailed study of this matrix-vector logic as a logic that implies algebraic operations between operators themselves, we refer to [32] and [35].

The matrix expressions for conjunction and disjunction immediately permit to define the matrices $S = NC$ and $P = ND$, corresponding to the Sheffer (or NAND) and the Peirce (or NOR) gates, respectively.

b3) **Implication.** The "material" implication corresponds in classical logic to the expression $p \to q \equiv \neg p \vee q$. The vector logic version of this equivalence leads to a matrix $L$ that represents vector logic "material" implication:

$$L = D(N \otimes I) \ .$$

As can be directly proved, the explicit expression for this implication is

$$L = s(s \otimes s)^T + n(s \otimes n)^T + s(n \otimes s)^T + s(n \otimes n)^T,$$

and the properties of classical implication are immediately verified:
$L(s \otimes s) = L(n \otimes s) = L(n \otimes n) = s$ and $L(s \otimes n) = n$.

b4) **The symmetric operators Equivalence and Exclusive-Or.**
In this vector logic the equivalence $p \equiv q$ corresponds to the following matrix:

$$E = s(s \otimes s)^T + n(s \otimes n)^T + n(n \otimes s)^T + s(n \otimes n)^T .$$

with, $E(s \otimes s) = E(n \otimes n) = s$ and $E(s \otimes n) = E(n \otimes s) = n$.

The Exclusive-Or is the negation of the equivalence, $\neg(p \equiv q)$; consequently it corresponds with the matrix $X = NE$ given by

$$X = n(s \otimes s)^T + s(s \otimes n)^T + s(n \otimes s)^T + n(n \otimes n)^T ;$$

hence $X(s \otimes s) = X(n \otimes n) = n$ and $X(s \otimes n) = X(n \otimes s) = s$ .

A simple numerical illustration of the form of these matrices for the case of $s = \begin{bmatrix} 1 & 0 \end{bmatrix}^T$ and $n = \begin{bmatrix} 0 & 1 \end{bmatrix}^T$ is included in the Appendix of reference [35]. This Appendix also illustrates the shape of logical matrices for 3-dimensional unit vectors that include a third vectorial truth-value.

## 2.2.2 Many-valued two-dimensional logic

In the case of two-valued vector logic, uncertainties in the truth-values can be introduced using vectors $f = \varepsilon s + \delta n$, with $\varepsilon, \delta \in [0,1]$, $\varepsilon + \delta = 1$. An interesting point is that these vectors can be directly processed by matrix operators that initially result from a "Boolean", non many-valued, logic. In this case, the many-valued character of the emerged logic has not been introduced *a priori* in the operators; instead, it is an *a posteriori* consequence of the uncertainties introduced in the inputs. The vectors of this many-valued logic map on scalar functions and generate a class probabilistic logic [30]. For vectors $u = \alpha s + \beta n$ and $v = \alpha' s + \beta' n$ the scalar many-valued (or "fuzzy") logic obtained from any two-valued matrix G is given by its projection over vector s :

$$\text{Val(scalars)} = s^T G(\text{vectors}),$$

with Val representing the scalar logical function associated with matrix G.

These projections produce the following results:
$$\text{NOT}(\alpha) = s^T N u = 1 - \alpha$$
$$\text{OR}(\alpha, \alpha') = s^T D (u \otimes v) = \alpha + \alpha' - \alpha\alpha'$$
$$\text{AND}(\alpha, \alpha') = s^T C (u \otimes v) = \alpha\alpha'$$
$$\text{IMPL}(\alpha, \alpha') = s^T L (u \otimes v) = 1 - \alpha(1 - \alpha')$$
$$\text{XOR}(\alpha, \alpha') = s^T X (u \otimes v) = \alpha + \alpha' - 2\alpha\alpha'$$

Using these equations we can define the corresponding negations:

$$\text{NOR}(\alpha, \alpha') = 1 - \text{OR}(\alpha, \alpha')$$
$$\text{NAND}(\alpha, \alpha') = 1 - \text{AND}(\alpha, \alpha')$$
$$\text{EQUI}(\alpha, \alpha') = 1 - \text{XOR}(\alpha, \alpha')$$

Let us define a set of probabilistic vectors

$$\Pi = \{ \gamma s + (1-\gamma) n : \gamma \in [0,1] \}.$$

An interesting point is that when the monadic or dyadic operators act over vectors belonging to this set, the output is also an element of this set. For monadic vectors, this result is obvious. Let us state the results for dyadic vectors in the following way:

*Lemma 2.1.*
If G is a dyadic logical matrix, and $u, v \in \Pi$, then $G(u \otimes v) \in \Pi$
*Proof.*
Let $u = \alpha s + (1-\alpha)n$, $v = \beta s + (1-\beta)n$, $u, v \in \Pi$. Now
$u \otimes v = \alpha\beta(s \otimes s) + \alpha(1-\beta)(s \otimes n) + (1-\alpha)\beta(n \otimes s) + (1-\alpha)(1-\beta)(n \otimes n)$
and $G(u \otimes v) = \phi_1 s + \phi_2 n$. But, necessarily, it must be
$\phi_1 + \phi_2 = \alpha\beta + \alpha(1-\beta) + (1-\alpha)\beta + (1-\alpha)(1-\beta) =$
$= [\alpha + (1-\alpha)][\beta + (1-\beta)] = 1$. Hence, $G(u \otimes v) = \phi_1 s + (1-\phi_1)n$, $\phi_1 \in [0,1]$ ∎

The fact that this set $\Pi$ is closed for any basic monadic or dyadic logical operation implies that it is also closed for any legal combination of matrix logical operations acting on probabilistic vectors and this produces a consistent many-valued logic. We are going to describe in the next Sections how this property allows extending the Boolean differential calculus to logic functions that go beyond the restrictions imposed by binary variables.

A wff of logical calculus can be always translated to the vector logic formalism provided that the logical variables (eg. 1,0 or True, False) are mapped on vectors s,n [30, 32]. We remark that s and n are vectors defined on $\mathbb{R}$ (or even over $\mathbb{C}$ [15]), and this fact determines the structure of the logical matrices. As was mentioned previously, once defined the matrix operators described in this Section some of the classical logical equivalences (equalities in the matrix-vector formalism) become intrinsic properties of the operators and do not depend on the structure of the input vectors [32]. Are examples of this fact the De Morgan Laws between conjunction and disjunction (with their duals between equivalence and exclusive-or) and the relation between implication and disjunction. Finally, let us comment that some basic tautologies (eg. *excluded middle*) are strictly valid only for vectors s and n. However, if the inputs of these tautologies are probabilistic vectors belonging to $\Pi$, the probabilistic weigh for s in the output is confined into the interval [(3/4), 1]. Let us denominate the outputs inside this interval *quasi-s* and its negation *quasi-n*.

3. BOOLEAN DERIVATIVES.
In what follows we will use a standard definition for the Boolean derivative [9, 18, 43]. Given a function

$y = f(x_1, \ldots, x_i, \ldots, x_p)$

with p Boolean variables $x_j \in \{0,1\}$, with $j=1,\ldots,p$, we can define a partial derivative of y respect to $x_i$ by means of the following equation:

$$\frac{\partial y}{\partial x_i} = \text{XOR}\left[f(x_1,\ldots,1,\ldots,x_p), f(x_1,\ldots,0,\ldots,x_p)\right].$$

This partial derivative is a propositional function. This fact allows to iterate the operation and, in this way, to obtain second and higher partial derivatives,

$$\frac{\partial^2 y}{\partial x_i \partial x_j}, \frac{\partial^3 y}{\partial x_i \partial x_j \partial x_k}, \text{ etc.}$$

Exclusive-or is used here as a discrete version of a differential operator. In the Boolean domain this operator indicates, in terms of truth-values, if there is any difference between the values of its arguments. Remark that this operator is symmetric, and it is not capable of attributing a sign to the variation.

## 4. THE DERIVATIVE OF MATRIX-VECTOR LOGICAL OPERATORS

In what follows, we represent a matrix logical operation of any complexity using the symbolical expression Op(u) and a vectorial variable $u \in \prod$. Hence, Op(u) describes expressions as diverse as $Mu, L(u \otimes v)$ and $L[C(u \otimes v) \otimes w]$.

*Definition 4.1. Boolean derivative of operator Op(u).*

$$\frac{\partial \text{Op}(u)}{\partial u} = X[\text{Op}(s) \otimes \text{Op}(n)]$$

This is, by definition, a partial derivative. If the logical expression involves other variables, they retain their own values. The operator X is the exclusive-or matrix previously described.

Remark that the Boolean derivative respect to variable u generates a logical function in which this variable disappears. Hence, the successive derivatives provoke an increasing reduction of complexity in the operations ("complexity" means, in this context, the number of logical variables)..

### 4.1. First Derivatives: A Case Study

The exclusive-or operator X defines a symmetrical function for all $u, v \in \prod$. The following properties are important in the evaluation of the Boolean derivatives:

$$X(u \otimes s) = Nu,$$
$$X(u \otimes n) = u.$$

(a) *Derivatives of monadic operations.*
Note that u = Iu. Hence,
$$\frac{\partial Iu}{\partial u} = \frac{\partial u}{\partial u} = X(s \otimes n) = s \ .$$
On the other hand, we have
$$\frac{\partial Nu}{\partial u} = X(n \otimes s) = s \ .$$
We now evaluate the derivatives of "constant" operations:

$$\frac{\partial Ku}{\partial u} = X(s \otimes s) = n \ ,$$
$$\frac{\partial Mu}{\partial u} = X(n \otimes n) = n \ .$$

Hence, there are some analogies with the classical derivatives of real analysis: $\partial Iu/\partial u = \partial Nu/\partial u = \partial u/\partial u = s$ is analogous to dx/dx = 1 and $\partial Ku/\partial u = \partial Mu/\partial u = n$ is equivalent to the classical result d(Constant)/dx = 0. The first situation illustrates the absence of signs in the logical derivatives.

(b) *Derivatives of dyadic operations.*
The above mentioned properties of operator X allows an immediate evaluation of the properties of the first derivatives of basic dyadic operations. In the following Table 3 we show these derivatives.

**Table 3 : First derivatives for matrix logical operators**

| Op | $\frac{\partial Op(u, v)}{\partial u}$ | $\frac{\partial Op(u, v)}{\partial v}$ |
|---|---|---|
| C | v | u |
| D | Nv | Nu |
| L | Nv | u |
| S | v | u |
| P | Nv | Nu |
| E | $X(v \otimes Nv)$ | $X(u \otimes Nu)$ |
| X | $X(v \otimes Nv)$ | $X(u \otimes Nu)$ |

We now state some remarks. The derivatives of these dyadic functions depend on the rest of the variables. For instance, in the case $\partial C/\partial u = v$, if v = s the transition of u from s to n provokes a modification of the logical value of conjunction C; instead, if v = n, the evaluation executed by C is unsensitive to the transition of u between s and n. The fact that the variables belong to the set $\Pi$ implies that these Boolean derivatives became capable of generating a fuzzy evaluation (even if the variable concerned in the differentiation transits a discrete step from s to n). Also note the asymmetry of the derivatives of the implication L. Finally, note that the derivatives of operators E and X generate quasi-tautologies. This can be easily seen by evaluating the scalar projections of one of these derivatives:

$$s^T \frac{\partial E}{\partial u} = s^T X(v \otimes Nv) = s^T E(v \otimes v) = \beta^2 + (1-\beta)^2 \equiv f(\beta) \ ;$$

$f(1) = f(0) = 1$ and $\min f(\beta) = \frac{1}{2}$ for $\beta = \frac{1}{2}$. This implies that for $\beta \in (0,1)$ the equivalence and the inequivalence are symmetrically sensitive to the modifications of their arguments, with the scalar projections of the outputs remaining into the upper half of the interval $(0,1)$.

(c) *Negation Lemmas.*
Given the operation Op(u), we state the following Lemma:

*Lemma 4.1*

$$\frac{\partial Op(u)}{\partial u} = \frac{\partial N Op(u)}{\partial u} \ , \ u \in \Pi$$

*Proof*
The equality is an immediate consequence of the identity $X = X(N \otimes N)$. ∎

This Lemma explains the coincidence of the derivatives of the pairs C and S, D and P, and E and X. It also explains why $\partial u/\partial u = \partial Nu/\partial u = s$. The second negation Lemma is a consequence of the symmetry of the operator X(u,v):

*Lemma 4.2*

$$\frac{\partial Op(u)}{\partial u} = \frac{\partial Op(Nu)}{\partial u} \ , \ u \in \Pi$$

*Proof*
It is immediate. ∎

A third negation Lemma can be stated as follows.

*Lemma 4.3*

$$\frac{\partial \text{Op}(u)}{\partial Nu} = \frac{\partial \text{Op}(u)}{\partial u} \ , \ u \in \Pi$$

*Proof*
$\partial \text{Op}(u)/\partial Nu = \partial \text{Op}(Nu')/\partial u'$ with $u' = Nu$. But
$\partial \text{Op}(Nu')/\partial u' = X[\text{Op}(Ns) \otimes \text{Op}(Nn)] = X[\text{Op}(n) \otimes \text{Op}(s)] =$
$= \partial \text{Op}(u)/\partial u$
due to the symmetry of X. ∎

It is interesting to note that Lemma 4.1 implies that

$$\frac{\partial N\text{Op}(u)}{\partial u} \neq N \frac{\partial \text{Op}(u)}{\partial u} \ .$$

Let us denominate "logical linearity" a situation in which an operation $\frac{\partial}{\partial u}$ satisfies:

(a) $\frac{\partial}{\partial u} C[t \otimes \text{Op}(u)] = C\left[t \otimes \frac{\partial}{\partial u} \text{Op}(u)\right]$, $t \in \{s, n\}$

(b) $\frac{\partial}{\partial u} X[\text{Op}(u) \otimes \text{Op}'(u)] = X\left[\frac{\partial}{\partial u} \text{Op}(u) \otimes \frac{\partial}{\partial u} \text{Op}'(u)\right]$ .

In this logical framework this pair of expressions is a version of the scalar linearity defined by $F(kx) = kF(x)$ and $F(x + y) = F(x) + F(y)$. It can be proved that the Boolean derivative satisfies property (a) and (b). The proof of property (a) holds immediately from the definition of the Boolean derivative. The proof of (b) requires proving the following matrix version of Abel bisymmetry equation:

$$X\big[X(a \otimes b) \otimes X(c \otimes d)\big] = X\big[X(a \otimes c) \otimes X(b \otimes d)\big]$$

for $\text{Op}(s) = a$, $\text{Op}'(s) = b$, $\text{Op}(n) = c$, and $\text{Op}'(n) = d$. A detailed proof is provided in Section 5, Lemma 5.2 in the context of cross derivatives.

The following rule for the logical product for variables $u, v, w, z \in \Pi$ (that can be directly proved from the properties of C and X) is similar to the rule of ordinary product differentiation in functions of real variables:

$$\frac{\partial}{\partial v} C\left[C(u \otimes v) \otimes w\right] = C(u \otimes w)$$

$$\frac{\partial}{\partial w} C\left\{C\left[C(u \otimes v) \otimes w\right] \otimes z\right\} = C\left[C(u \otimes v) \otimes z\right].$$

Chain rule is not fully satisfied in the case of Vichniac formalism, except along a path that is generated during the time-evolution of the Boolean elementary cellular automata and that he denominates "light cone" [43]. The conditions of validity of the chain rule inside our matrix-vector formalism require further investigation. Here we can comment that there are expressions that verify the chain rule under different conditions of generality. We illustrate this point with two examples:

*Example 1.* $F = E\left[L(u \otimes v) \otimes X(w \otimes w)\right]$

The derivative respect to v is

$$\frac{\partial F}{\partial v} = X\left\{X(w \otimes w) \otimes E\left[Nu \otimes X(w \otimes w)\right]\right\}.$$

Applying the chain rule, we obtain

$$C\left(\frac{\partial F}{\partial L} \otimes \frac{\partial L}{\partial v}\right) =$$

$$C\left\{E\left[X(w \otimes w) \otimes X(w \otimes w)\right] \otimes u\right\}.$$

These expressions are not equivalent in general, but if $w \in \{s, n\}$ they satisfy the equality

$$\frac{\partial F}{\partial v} = C\left(\frac{\partial F}{\partial L} \otimes \frac{\partial L}{\partial v}\right) = u$$

*Example 2.* $F = L\left[C(u \otimes v) \otimes D(w \otimes w)\right]$

Now the derivative respect to v is

$$\frac{\partial F}{\partial v} = NL\left[u \otimes D(w \otimes w)\right] = NL(I \otimes D)(u \otimes w \otimes w).$$

Applying the chain rule, we obtain

$$C\left(\frac{\partial F}{\partial L} \otimes \frac{\partial L}{\partial v}\right) = C\left[u \otimes ND(w \otimes w)\right] = C(I \otimes ND)(u \otimes w \otimes w)$$

We used the Kronecker product factorizations in both expressions. It can be easily proved that $NL(I \otimes D) = C(I \otimes ND)$ and, consequently, in this Example 2 the chain rule is valid for probabilistic vectors $u, v, w \in \Pi$.

Finally, we pay attention to the Leibniz rule. In his formalism for pure Boolean functions, Vichniac showed that Leibniz rule exhibit an additional correction term. In the case of our vector logic formalism, an elementary Leibniz rule with its classic structure, ie, $d[f(x).g(x)]/dx = [df(x)/dx].g(x) + f(x).[dg(x)/dx]$, cannot exists because the logical variable vanishes during differentiation. Potential modifications of the differentiation procedures in order to re-obtain the Leibniz rule deserve further explorations

5. CROSS DERIVATIVES
It is a remarkable fact that, even though it is not a complete linear operation, the Boolean derivative retains a certain parallelism with the partial derivatives of classical differential calculus. In particular, as we show in what follows, the equality between cross derivatives holds. To prove this equality we need previously to prove two Lemmas:

*Lemma 5.1*
Given the logical operation $G(u,v)$, with $u, v \in \Pi$ and G being a legal (but otherwise arbitrary) logical expression built up using the operators of vector logic, to prove the equality

$$\frac{\partial}{\partial v}\left(\frac{\partial G(u,v)}{\partial u}\right) = \frac{\partial}{\partial u}\left(\frac{\partial G(u,v)}{\partial v}\right)$$

requires to prove the equality

$$X\left[X(a \otimes b) \otimes X(c \otimes d)\right] = X\left[X(a \otimes c) \otimes X(b \otimes d)\right]$$

with $a = G(s,s)$, $b = G(n,s)$, $c = G(s,n)$ and $d = G(n,n)$.

Note that a, b, c, or d can be probabilistic outputs if $G(u,v)$ represents a function involving other probabilistic variables, eg. $G(u, v) = H(u, v, w) = L\left[NC(u \otimes v) \otimes w\right]$.

*Proof*
It is immediate, taking into account that

$$\frac{\partial}{\partial v}\left(\frac{\partial G(u,v)}{\partial u}\right) = \frac{\partial}{\partial v}X\left[G(s,v) \otimes G(n,v)\right] =$$

$$= X\{X[G(s,s) \otimes G(n,s)] \otimes X[G(s,n) \otimes G(n,n)]\}$$

and

$$\frac{\partial}{\partial u}\left(\frac{\partial G(u,v)}{\partial v}\right) = \frac{\partial}{\partial u} X[G(u,s) \otimes G(u,n)] =$$

$$= X\{X[G(s,s) \otimes G(s,n)] \otimes X[G(n,s) \otimes G(n,n)]\}. \quad \blacksquare$$

We now establish the following Lemma:

*Lemma 5.2*
If $a, b, c, d \in \Pi$ then

$$X\big[X(a \otimes b) \otimes X(c \otimes d)\big] = X\big[X(a \otimes c) \otimes X(b \otimes d)\big]$$

*Proof*
The outputs of the last applications of the operator X have the form $\varphi s + (1-\varphi)n$, with $\varphi \in [0,1]$, because its inputs also have this form, as is assured by Lemma 2.1. Hence, if we have

$$a = \alpha s + (1-\alpha)n \quad , \quad b = \beta s + (1-\beta)n \;,$$
$$c = \gamma s + (1-\gamma)n \quad , \quad d = \delta s + (1-\delta)n \;,$$

the coefficient $\varphi$ satisfies the following functional equation:

$$\varphi[\varphi(\alpha,\beta), \varphi(\gamma,\delta)] = \varphi[\varphi(\alpha,\gamma), \varphi(\beta,\delta)] \;.$$

This functional equation is a particular case of the "bisymmetry equation" [1]. For the case of exclusive-or operator X, the function $\varphi$ has the form $\varphi(x,y) = x + y - 2xy$. Using this expression, and developing both members of the precedent bisymmetry equation, the identity can be directly proved. $\blacksquare$

We can now establish the equality between cross derivatives:

*Theorem 5.1*
Given a logical function G(u,v) with $u, v \in \Pi$, then

$$\frac{\partial}{\partial v}\left(\frac{\partial G(u,v)}{\partial u}\right) = \frac{\partial}{\partial u}\left(\frac{\partial G(u,v)}{\partial v}\right)$$

*Proof*
From Lemma 2.1 we know that $G(u,v) \in \prod$. In addition, we are under the conditions of the hypothesis of Lemma 5.1; hence, following Lemma 5.2 the equality holds. ∎

We will use for these equivalent derivatives the notation $\partial^2 G(u,v)/\partial[u,v]$. In the following table we show the cross derivatives of the basic operators.

**Table 4 : Cross derivatives for matrix logical operators**

| Op | $\partial^2 Op(u,v)/\partial[u,v]$ |
|---|---|
| C | s |
| D | s |
| L | s |
| S | s |
| P | s |
| E | n |
| X | n |

## 6. SUCCESSIVE DERIVATIVES

After the first derivative, the concerned variable disappears. Due to this fact, one can assume that the successive derivatives respect to this variable are irrelevant. We mention that, in fact, this is the case when the Boolean derivatives act over binary variables (for a recent reference see Cheng, proposition 3.2 [9]), Nevertheless, we are going to explore in this section a heuristic definition of the successive Boolean derivatives that shows that the logical functions in the vectorial domain can display some unexpected properties.

Given a logical operator that depends on three vectors $Op_1(u,v,w)$, with
$u = \alpha s + (1-\alpha)n$, $v = \beta s + (1-\beta)n$, $w = \gamma s + (1-\gamma)n$ we have:

$$\frac{\partial Op_1(u,v,w)}{\partial u} = Op_2(v,w) = f(\beta,\gamma)s + [1-f(\beta,\gamma)]n \ .$$

Here $Op_1(u,v,w)$ represents a logical function of three variables defined using the basic dyadic or monadic operators (eg. $L[X(u \otimes v) \otimes w]$ or $E[C(u \otimes v) \otimes D(v \otimes Nw)]$). We only analyze the three-variable situation, since the generalization is immediate.

We define the second derivative as follows:

$$\frac{\partial}{\partial u}\frac{\partial Op_1(u,v,w)}{\partial u} = \frac{\partial^2 Op_1(u,v,w)}{\partial u^2} =$$
$$X[Op_2(v,w) \otimes Op_2(v,w)] = f'(\beta,\gamma)s + [1-f'(\beta,\gamma)]n,$$

with

$$f'(\beta,\gamma) = 2f(\beta,\gamma)[1-f(\beta,\gamma)].$$

This last expression is based on the following Lemma:

*Lemma 6.1*
For $z, u \in \Pi$, with z being independent of u, the derivative of z respect to u exists, and is given by

$$\frac{\partial z}{\partial u} = X(z \otimes z)$$

*Proof*
Note that for all $z, u \in \Pi$ we have
a)  $Ku = s$
b)  $C(s \otimes z) = z$.

Consequently, any expression that is not a function of u, can be converted into a "silent" function of u of the form $z = C(Ku \otimes z)$ that allows to apply the algorithm of differentiation with respect to u; in our case

$$\frac{\partial z}{\partial u} = \frac{\partial C(Ku \otimes z)}{\partial u} =$$
$$X[C(Ks \otimes z) \otimes C(Kn \otimes z)] = X(z \otimes z) \quad \blacksquare$$

In the Boolean domain, being z equal to s or n, the second derivative gives n, the expected result for the derivative of a Boolean constant function. This fact supports the agreement of the result of Lemma 6.1 with our expectations for the binary domain. But if the other variables of the logical functions are out of the Boolean domain, the fuzziness imposes an increase of uncertainty. Here is the argument: If a first derivative with respect to u has the structure $\text{Op}_2(v,w) = f(\beta,\gamma)s + [1-f(\beta,\gamma)]n$ then, after a second differentiation with respect to u, for any $f(\beta,\gamma)$ we obtain the following mapping:

$$f'(\beta,\gamma) = 2f(\beta,\gamma)[1-f(\beta,\gamma)] .$$

A mapping of this type, with the general form $\varepsilon' = 2\varepsilon(1-\varepsilon)$, has two fixed points, $\varepsilon = 0$ and $\varepsilon = 1/2$. The successive applications of this mapping define a dynamical system in which the point $\varepsilon = 0$ is unstable, and $\varepsilon = 1/2$ is a global attractor. Hence, if the function $f(\beta,\gamma)$ exhibits some degree of fuzziness ($f(\beta,\gamma) \neq 0,1$) after the first derivative, the successive higher derivatives "push" the system towards $f^{(n \to \infty)}(\beta,\gamma) = 1/2$, the value of maximum uncertainty.

## 7. TWO APPLICATIONS OF THE LOGICAL DERIVATIVES

In what follows we describe two situations that allow to illustrate how the Boolean derivative, applied over matrix logical functions, can produce some interesting results.

### 7.1. Simplification of Logical Chains

We denominate "logical chain" any well-formed logical formula involving more than one logical operator. In the following example we analyze the way in which the Boolean derivative provides us with explicit vectorial expressions. These expressions are useful to analyze the sensitivity of the logical chain to its different variables.

*Example:* $(p \vee q) \to (\neg q \wedge p)$

Representing the truth-values of propositions p and q by the probabilistic vectors u and v, respectively, we have the following vectorial representation:

$$h = L\big[D(u \otimes v) \otimes C(Nv \otimes u)\big] .$$

The derivatives are

$$\frac{\partial h}{\partial u} = X(v \otimes v),$$
$$\frac{\partial h}{\partial v} = L(u \otimes u).$$

Notice that in the Boolean domain the derivative with respect to u, $X(u \otimes u)$, gives the negation vector n, and in the probabilistic domain it provides vector biased towards n, with $s^T \partial h/\partial u \leq 0.5$. The meaning of this result is that the "argument" represented by the logical chain has low sensitivity to the truth-value of the variable u. On the contrary, it is sensitive to v: in the case in which v = s, the value of h is n; when v = n, the value of h is $L(u \otimes u)$ (this is a version of the *excluded middle* that gives s for u Boolean and a quasi-s for a probabilistic u). Consequently, the chain represented by h practically collapses in the following approximate expression:

$$h \approx Nv$$

(for the Boolean values u = s and u = n, the expression becomes a true equality).
In words: "if he is a good person or he is smart, then, he is not smart and he is a good person" is an obscure argument approximately equivalent to the proposition "he is not smart". This proposition can be true or false in the case the truth-value of the variable v is Boolean, with a weighted uncertainty if v is probabilistic.

7.2. **The Derivatives of Basic Tautologies**
In the framework of the classical logic, the tautologies are logical expressions that produce the truth-value "true" for all the possible values of their logical variables. Some of the most basic tautologies are the following:
*Excluded middle*: $p \vee \neg p$

*Modus ponens*: $[p \wedge (p \to q)] \to q$

*Hypothetical syllogism*: $[(p \to q) \wedge (q \to r)] \to (p \to r)$

In the formalism of vector logic, these tautologies can be respectively expressed as matrix-vector operations in the following way:

$$EM(u) = D(u \otimes Nu)$$
$$MP(u, v) = L\{C[u \otimes L(u \otimes v)] \otimes v\}$$
$$HS(u, v, w) = L\{C[L(u \otimes v) \otimes L(v \otimes w)] \otimes L(u \otimes w)\}$$

When u,v,w are probabilistic vectors, the scalar projection of these expressions remains bounded to the interval $[(3/4), 1]$.

The application of the logical derivative to these matrix expressions generates an interesting result: derivatives of *hypothetical syllogism* produce the negation of *modus ponens*, and derivatives of *modus ponens* produce the negation of the *excluded middle*. In what follows we calculate the derivatives of the excluded middle and the *modus ponens,* and we show the results for the case of *hypothetical syllogism.*

*Case 1. Excluded middle*
The vectorial version of the *excluded middle* implies the following equality:

$$\mathrm{EM}(u) = D(u \otimes Nu) = D(Nu \otimes u) = L(u \otimes u).$$

The derivative is given by

$$\frac{\partial \mathrm{EM}(u)}{\partial u} = X\left[L(s \otimes s) \otimes L(n \otimes n)\right] = n = Ns.$$

*Case 2. Modus ponens*
In this case we have

$$\frac{\partial \mathrm{MP}(u,v)}{\partial u} = X\left\{L\left\{C\left[s \otimes L(s \otimes v)\right] \otimes v\right\} \otimes L\left\{C\left[n \otimes L(n \otimes v)\right] \otimes v\right\}\right\} =$$
$$X\left[L(v \otimes v) \otimes L(n \otimes v)\right] = X\left[L(v \otimes v) \otimes s\right] = NL(v \otimes v).$$

Hence,

$$\frac{\partial \mathrm{MP}(u,v)}{\partial u} = ND(Nv \otimes v) = N\left[\mathrm{EM}(v)\right].$$

We also have

$$\frac{\partial \mathrm{MP}(u,v)}{\partial v} = ND(u \otimes Nu) = N\left[\mathrm{EM}(u)\right].$$

The cross derivatives are

$$\frac{\partial^2 \mathrm{MP}(u,v)}{\partial [u,v]} = n.$$

*Case 3. Hypothetical syllogism*
For the *hypothetical syllogism,* it can be proved that

$$\frac{\partial \, HS(u,v,w)}{\partial \, u} = N\big[MP(v,w)\big] \, ,$$

$$\frac{\partial \, HS(u,v,w)}{\partial \, w} = N\big[MP(Nu, Nv)\big].$$

The derivative respect to the "pivot" vector v produces a complex expression:

$$\frac{\partial \, HS(u,v,w)}{\partial \, v} = X\big\{L\big[w \otimes L(u \otimes w)\big] \otimes L\big[Nu \otimes L(Nw \otimes Nu)\big]\big\}.$$

Considering that X = NE, the expression

$$E\big\{L\big[w \otimes L(u \otimes w)\big] \otimes L\big[Nu \otimes L(Nw \otimes Nu)\big]\big\}$$

defines a tautology in the binary domain, as can be easily proved evaluating this expression for $u, v \in \{s, n\}$. From the identity $L(I \otimes L) = L(C \otimes I)$ (see [32]), it follows that $L\big[a \otimes L(a \otimes b)\big] = = L\big[C(a \otimes b) \otimes a\big]$, with $a, b \in \prod$. This is the version, in the vector logic representation, of other of the classical tautologies: $(p \wedge q) \to p$ that corresponds to the matrix expression $L\big[C(a \otimes b) \otimes a\big]$. Likewise, it is clear that within the dyadic calculus if A and B are tautologies then $A \equiv B$ is also a tautology. Consequently the derivative of the *hypothetical syllogism* with respect to w generates the negation of a new tautology TD. This new tautology TD results from the equivalence of two "subtautologies" of the form $(p \wedge q) \to p$.

An interesting point is that the derivatives of the subtautology

$$ST(u,w) = L\big[C(w \otimes u) \otimes w\big]$$

are:

$$\frac{\partial \, ST(u,w)}{\partial \, u} = NL(w \otimes w)$$

$$\frac{\partial \, ST(u,w)}{\partial \, w} = n$$

Finally, we conclude that the pivot derivative of the *hypothetical syllogism* can be

represented by

$$\frac{\partial\, HS(u, v, w)}{\partial v} = F(u, w)$$

with

$$F(s, w) = NL(w \otimes w) = N[EM(w)],$$
$$F(u, s) = n = Ns,$$
$$F(n, w) = n = Ns,$$
$$F(u, n) = N[EM(u)].$$

In conclusion, the results presented in this section show us that the derivatives of some basic tautologies generate tautologies of inferior order (dimension) pre-multiplied by the negation matrix. We can represent this decreasing hierarchy in the following transition diagram:

$$HS(u, v, w) \xrightarrow{N\partial/\partial u} MP(v, w) \xrightarrow{N\partial/\partial v} EM(w) \xrightarrow{N\partial/\partial w} s$$

The lowest order tautology is vector s. The results obtained in this Section illustrate clearly how when the Boolean derivatives are applied to semantically meaningful logical operators they can generate packed expressions susceptible of the same kind of semantic logical interpretation. Obviously, each one of the expressions shown in this section about tautologies can be developed according with the syntax of typical Boolean derivatives, losing their packed organization.

## 8. AN INTEGRAL CALCULUS FOR LOGICAL OPERATORS

In this Section we are going to describe a class of logical function that presents a formal analogy with the indefinite integral of real analysis.

*Definition 8.1. Boolean Integral.*
Given a logical vectorial function $Op$, we define its Boolean integral as another logical function $\Upsilon$ such that

$$\frac{\partial \Upsilon}{\partial \tau} = Op,$$

where $\tau$ is a new vectorial variable not included in $Op$. We use the following notation:

$$\Upsilon = \int \mathrm{Op}\, \partial \tau$$

(we define $\partial \tau$ as a Boolean differential).

This Boolean integral extends the dimensionality of the domain of the operation Op (the contrary of Boolean derivation, that eliminates the variable with respect to which derivation is made). As we show in the next Theorem, it is always possible to associate a family of Boolean integrals to any arbitrary logical function Op.

*Theorem 8.1. General Integral.*
An arbitrary logical operation Op admits a general Boolean integral of the form

$$\int \mathrm{Op}\, \partial \tau = \mathrm{HL}(\mathrm{Op} \otimes \mathrm{H}'\tau) \;,\quad \mathrm{H}, \mathrm{H}' \in \{\mathrm{I}, \mathrm{N}\}\;,$$

$\tau$ being a logic vector.

*Proof*
It is immediate using the Definition 8.1 and considering that $\mathrm{L}(\mathrm{Op} \otimes \mathrm{s}) = \mathrm{s}$, $\mathrm{L}(\mathrm{Op} \otimes \mathrm{n}) = \mathrm{NOp}$ and $\mathrm{X}(\mathrm{NOp} \otimes \mathrm{s}) = \mathrm{X}(\mathrm{s} \otimes \mathrm{NOp}) = \mathrm{NNOp} = \mathrm{Op}$. If H = N we are under Lemma 4.1. ∎

This general integral admits for different values of H and H' the following versions:

1) $\int \mathrm{Op}\, \partial \tau = \mathrm{L}(\mathrm{Op} \otimes \tau)$

2) $\int \mathrm{Op}\, \partial \tau = \mathrm{NL}(\mathrm{Op} \otimes \tau)$

3) $\int \mathrm{Op}\, \partial \tau = \mathrm{NL}(\mathrm{Op} \otimes \mathrm{N}\tau) = \mathrm{C}(\mathrm{Op} \otimes \tau)$

4) $\int \mathrm{Op}\, \partial \tau = \mathrm{L}(\mathrm{Op} \otimes \mathrm{N}\tau) = \mathrm{NC}(\mathrm{Op} \otimes \tau)$

*Corollary 8.1.*

$$\int \mathrm{NOp}\, \partial \tau = \mathrm{L}(\mathrm{NOp} \otimes \mathrm{N}\tau) = \mathrm{L}(\tau \otimes \mathrm{Op})$$

Contraposition of implication L is valid for probabilistic vectors, as can be directly proved given the identity $L = D(N \otimes I)$ and the commutativity of disjunction D:

$$L(u \otimes v) = D(N \otimes I)(u \otimes v) = D(Nu \otimes v) = D(v \otimes Nu) =$$
$$D(NNv \otimes Nu) = D(N \otimes I)(Nv \otimes Nu) = L(Nv \otimes Nu), \quad u, v \in \Pi .$$

This Corollary can be used to show a kind of non-linearity of this integral respect to N because pre-multiplying by N each one of the previous equalities 1) to 4) we can verify that

$$\int NOp(u)\partial \tau \neq N \int Op(u)\partial \tau .$$

The general integral is independent of the particular form of the function Op.

Apart from this integral, it is possible to define particular integrals directly dependent on the form of Op. In what follows, to organize the argument, we label the position of a variable into a given logical expression with a number. Using these labels, we can define the substitutions that create the particular integrals. Let $Op(u[1], u[2], \ldots, u[n])$ be a logical function where the variable u[i] fills the position i in the structure of the logical expression Op. We remark that u[i] and u[j], $i \neq j$, can be the same variable (eg, for $C[L(v \otimes w) \otimes Nv]$ we have u[1] = v, u[2] = w, u[3] = v). A heuristic procedure to obtain a particular integral consists in the substitution

$$u[i] \rightarrow B_i(u[i], \tau) ,$$

where $B_i(u[i], \tau)$ is a logical operator. The idea is to look for transformations of variables able to force a detachment during differentiation. We base this procedure in the following theorem:

*Theorem 8.2.*
The substitutions $u \rightarrow B(u, \tau)$ and $v \rightarrow B'(v, \tau)$ inside a logical function F(u,v) generate a particular integral of such function in the following cases:

(c1) $B(u,s) = u$ ; $B'(v,s) = v$ ; $F[B(u,n), B'(v,n)] = n$ ,

(c2) $B(u,n) = u$ ; $B'(v,n) = v$ ; $F[B(u,s), B'(v,s)] = n$ .

*Proof*
The evaluation of the Boolean derivative with the substituted variables show that both (c1) and (c2) produce the detachment condition

$$X[F(u,v) \otimes n] = X[n \otimes F(u,v)] = F(u,v) ,$$

that assures that the enlarged function is a particular integral. ∎

For instance, the substitutions $u \to B(u,\tau)$ and $v \to B'(v,\tau)$ such that $B(u,s) = u$, $B(u,n) = n$ and $B'(v,s) = v$, $B'(v,n) = n$, generate particular integrals in the case where $Op = L(u \otimes v)$ (see the following Example 3).

In the following we give three examples to illustrate the procedure. In each example we describe for the same operation the general and the particular integrals, indicating the particular integrals by the subscript P.

*Example 1*   $Op = D(u \otimes Nu)$
a) General Integral:

$$\Upsilon = \int D(u \otimes Nu) \partial \tau = L\big[D(u \otimes Nu) \otimes \tau\big]$$

b) Particular integral:
$u \to C(\tau \otimes u)$, $Nu \to C(\tau \otimes Nu)$

$$\Upsilon_P = \int_P D(u \otimes Nu) \partial \tau = D\big[C(\tau \otimes u) \otimes C(\tau \otimes Nu)\big]$$

The proof is immediate, because the derivative of this last expression is

$$\frac{\partial \Upsilon_P}{\partial \tau} = D(u \otimes Nu)$$

*Example 2*   $Op = L(u \otimes Nv)$
a) General Integral:

$$\Upsilon = L\big[L(u \otimes Nv) \otimes \tau\big]$$

b) Particular integral:
$$u \to D(u \otimes \tau), \quad v \to D(v \otimes \tau)$$
Hence,

$$\Upsilon_P = L\big[D(u \otimes \tau) \otimes ND(v \otimes \tau)\big].$$

*Example 3*   $Op = L(u \otimes v)$

a) General Integral:

$$\Upsilon = L\big[L(u \otimes v) \otimes \tau\big]$$

b) Particular integral:
$$u \to L(\tau \otimes u), \quad v \to C(\tau \otimes v)$$
Hence,

$$\Upsilon_P = L\big[L(\tau \otimes u) \otimes C(\tau \otimes v)\big].$$

It is important to emphasize that different substitutions can lead to different particular integrals. In the case of Example 3, if we perform the substitution

$$u \to E\big[\tau \otimes C(\tau \otimes u)\big], \quad v \to C\big[\tau \otimes E(\tau \otimes v)\big]$$

the following expression is also a particular integral for the function $L(u \otimes v)$:

$$\Upsilon'_P = L\big\{E\big[\tau \otimes C(\tau \otimes u)\big] \otimes C\big[\tau \otimes E(\tau \otimes v)\big]\big\}.$$

## 9. PERSPECTIVES

The Theorem 8.2 described in the previous Section illustrates an interesting case of "conditioned heuristics", where the invention of the functions $B(u,\tau)$ and $B'(v,\tau)$ is modulated by rigid impositions. In fact, this kind of invention constrained by rigid conditions is also present in the classical deduction systems, that operate using substitution and detachment. In that sense, let us mention the important observations by Łukasiewicz in his 1931 article on "generalizing deduction" [26]. There, he describes how deduction can be concomitant with an increase in the generality (or the complexity) of the expressions. In [26, p.191] he wrote : *"[...] we have demonstrated that in certain cases we can pass, in a deductive manner, from the particular to the general [...]"*. The suggestive discover of Łukasiewicz and this comment induces to think that integration and deduction can be related, a point that deserves further investigation.

As we saw in Section 2, the vectorial truth-values can be combined to produce a representation of uncertainties by mean of probabilistic weights. The informational aspect of this calculus becomes clear. On the one hand, the loss of information during differentiation is evident due to the reduction of the variables. In the limit of the successive differentiation of a logical vector, as was shown in Lemma 6.1, this calculus produces an effect interpretable as an increase of the uncertainties of the remaining logical variables provided that these variables are probabilistic vectors. On the other hand, the integral calculus shows an increase in the number of variables and in the case of particular integrals, a potential complexification of the logical functions.

We can ask if this calculus of logical operators could have some practical importance aside their potential theoretical interest. In this sense, it is interesting to note that the first and the second derivatives described in Tables 3 and 4, separates clearly the exclusive-or and the equivalence from the other basic dyadic functions. This fact is remarkable because in the case of elementary cellular automata this two operations are the main responsible of the generation of dynamic complexity [33, 34]. Hence this calculus can be a way to penetrate the hard problems of the generation of complexity in formal systems and in dynamic models based in logical operations. In addition, in a similar way as happens in switching exploration of circuits' behaviour, this calculus can be useful to explore the sensitivity of a complex logical reasoning to the different arguments involved, and in this way evaluate the relevance or the irrelevance of parts of the argument, as in our miniature example of Section 7.1.

**Aknowledgments**
The author acknowledges the partial financial support by PEDECIBA and CSIC-UdelaR. He wishes to thank both reviewers for helpful comments and useful insights.